# Engagement Zones for a Turn Constrained Pursuer

Thomas Chapman[1], Isaac E. Weintraub[1], Alexander Von Moll[1], and Eloy Garcia[1]

*Abstract*—This work derives two basic engagement zone models, describing regions of potential risk or capture for a mobile vehicle by a pursuer. The pursuer is modeled having turn-constraints rather than simple motion. Turn-only (C-Paths) and turn-straight (CS-Paths) paths are considered for the pursuer of limited range. Following the derivation, a simulation of a vehicle avoiding the pursuer's engagement zone is provided.

## I. INTRODUCTION

In path planning, the objective is to obtain a trajectory through space that optimizes some performance criterion subject to a set of constraints. This is a relevant problem in flying autonomous uncrewed aerial systems as they rely upon spatio-temporal information to move about without direct human interference. The task of path planning in congested airspace is a challenge and priority for enabling air mobility [1], [2]. Establishing regions of space where vehicle conflict could arise can aid in safe mission planning. This work considers the reachable range of a vehicle moving with constant speed and constrained turn radius and extends those regions to *engagement zones* (EZ) that account for the time of travel between a target vehicle and the other mobile agent.

In service of developing numerical tools to analyze optimal navigation given adversarial constraints, previous work [3] defined the concept of a *engagement zone* (EZ) - a dynamic region of space wherein an agent can be captured by a pursuer if it does not deviate from its current trajectory. This region assumes an adversarial relationship, which is a conservative approach, as inter-vehicle planning seeks to minimize interception.

Prior work considered the pursuer to have simple motion [3]; this work extends prior art to enforce a heading rate or minimum turn radius on the pursuer, more common in fixed-wing flight.

### A. Dubins Cars and Curvature Contrained Motion

The canonical model for this form of turn-constrained motion is the Dubins' car after the seminal work in [4], which has since been studied extensively [5], [6], [7], [8], both in the context of pursuit/evasion and of vehicular path planning.

Notably, the work examining pursuit-evasion involving a Dubins vehicle has been largely restricted to questions of equilibrium states of capture/escape - can the pursuer eventually intercept the target, or can the target indefinitely evade the pursuer given a set of initial conditions. For instance, [8] derived the conditions under which a pursuer could capture an evader where both moved with differing speeds and turn constraints in a plane. This work was extended to 3 dimensions in [6]. Little work has been dedicated to incorporating separate objectives for the target beyond simply avoiding interception.

### B. Reachability Analysis and Control

A *reachability region* is defined as a region of space wherein a performance limited vehicle can reach or be effective. In [9], the reachable sets of a particle moving with constant speed subject to curvature constraints were considered. In their work, the region was presented; but, no derivation was provided – in this work we provide such a derivation. In [10], the reachable region with bounded curvature is described; but offered an incomplete description of the reachability region. Other works considered reachable sets for car models [11], [12], and those with asymmetric control [13]. In more recent work, the use of machine learning was used to predict reachable sets [14].

The interception of a target by a Dubins' vehicle was examined in [15], which provided easy-to-compute lower and upper bounds on the minimum time to interception of a non-maneuvering target.

### C. Path Planning around Engagement Zones

The objective of planning motion around a single or field of EZs has been considered in [16], [17], [18]. In these works, cardioid models were used for lack of a model tied to first principles. The use of nonlinear programming techniques for driving and around and through cardioid EZs was considered in [16], planning around two EZs using pseudo-spectral methods was performed in [17], and sampling based approaches for avoiding a field of EZs was conducted in [18].

### D. Differential Games and Pursuit-Evasion

The consideration of adversarial engagements is a conservative one that offers benefits from a worst-case perspective. In this vein, games of pursuit-evasion and differential games [19] are of great value.

*This paper is funded in part by AFOSR, LRIR 24RQCOR002. DISTRIBUTION STATEMENT A. Approved for public release. Distribution is unlimited. AFRL-2025-0512; Cleared 29-Jan-2025.

[1]Chapman, Weintraub, Von Moll, Garcia are with Aerospace Systems Directorate, Wright-Patterson AFB, OH 45433, {`thomas.chapman.6`, `isaac.weintraub.1`, `alexander.von_moll`, `eloy.garcia.2`} `@us.af.mil`.

Two related classical differential games are "The game of Two Cars" and the "Homicidal Chauffeur." In the Game of Two Cars, two turn-constrained vehicles are pitted against one another; the aim is to evade or ensure capture in the adversarial engagement [20], [21], [22]. In the Homicidal Chauffeur differential game, a pedestrian aims to evade capture by a faster and turn-constrained pursuer [19].

*E. Contributions*

This work formulates two additional "basic" engagement zones (BEZ) - basic in that they derive from first principles rather than real system simulated capabilities - and compares them with prior art. We considers a turn-constrained pursuer. The pursuer is modelled as moving at a constant speed $v$ with a minimum turn radius of $\bar{a}$, and a maximum effective range (equivalently, time of flight $t$) given by $R = vt$. For simplicity, the pursuer's capture radius is taken to be 0, though a nonzero capture radius only serves to radially increase the EZs described herein by its value.

Two potential interception strategies are considered, one where the pursuer moves in a constant curve for the entirety of its motion, and one where it begins in a constant curve of some length, before breaking off and continuing in a straight line tangential to its previous circular motion.

Lastly, an optimal navigation problem is posed and solved, demonstrating the ability of EZ-informed motion to outperform more conservative approaches without compromising target safety.

## II. REACHABILITY REGIONS

In this section, two reachability regions of a pursuer are provided under the assumption that the pursuer has limited range. One strategy assumes the pursuer takes a single, constant turn referred to as a C-path. The second strategy has the pursuer moving along a constant turn and then a straight line segment, referred to as a CS-path. In both cases the turn's radii are assumed to be greater than or equal to a minimum $\bar{a}$.

*A. Reachability Region of C-Paths*

Consider that the vehicle only takes a single smooth turning trajectory. This could be the case if the vehicle only has the control authority to make one turn rather than having the ability to turn and then straighten out.

Let $a$ be the chosen turn radius of the vehicle such that $0 < \bar{a} \leq a$. The turn radius relates to the angle $\theta$ through which the vehicle turns:

$$a = \frac{vt}{\theta} \quad (1)$$

From Fig. 2 then, assuming that the vehicle initially points along the x-axis, the locus of possible terminal positions of the vehicle at maximum range is given by:

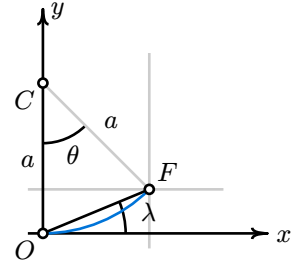

Fig. 2. A vehicle initially pointing along the $x$-axis takes a constant turn until reaching a point $F$.

$$x(\theta) = \frac{vt}{\theta} \sin(\theta)$$
$$y(\theta) = \frac{vt}{\theta}(1 - \cos(\theta)), \quad -vt/\bar{a} \leq \theta \leq vt/\bar{a} \quad (2)$$

Importantly, this is *not* polar $\theta$. For the purposes of this paper, the polar angle is denoted as $\lambda$.

There is also a convenient polar representation of (2):

$$r(\lambda) = \frac{vt}{\lambda}\sin(\lambda), \quad -\frac{vt}{2\bar{a}} \leq \lambda \leq \frac{vt}{2\bar{a}} \quad (3)$$

This shape, plotted in Fig. 3, is a cochleoid and can be derived from (2) and Fig. 2 with the observation that $\lambda = \frac{\theta}{2}$. The reachability region is all of the space contained radially between this bounding curve and the minimum turn radius circles within. The blue line shows a single trajectory for a smooth turning vehicle.

Note that the circles have been added separately, as they are not captured in (2) or (3), since those equations only describe the *endpoints* of all C-paths of length $vt$ with turn radius $a \geq \bar{a}$, not the full trajectories.

Also note that the locus of endpoints go inside the outer boundary when $vt$ exceeds $2\pi\bar{a}$. While accurate, this can interfere with the use of the boundary to define the reachability region, but this can be resolved by restricting the polar equation's domain in (3) to be $\lambda \in [-\pi, \pi]$ for large effective ranges.

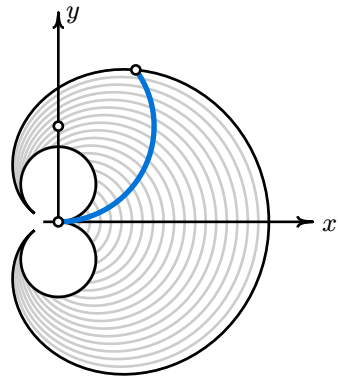

Fig. 3. A vehicle (blue) takes a smooth turn for a fixed amount of time. The locus of all possible terminal conditions that do not exceed the turn radius or range of the vehicle are shown in black. In grey are all the iso-cost lines for other time intervals to show how the C-Frontier grows with time.

## B. Reachability Region of CS-Paths

In the case where the vehicle takes a CS-path, the locus of endpoints of turn-straight curves of length $vt$ with minimum turn radius $\bar{a}$ may be obtained and are as described in [9] to be:

$$x(\theta) = \bar{a}\sin\theta + (vt - \bar{a}\theta)\cos\theta$$
$$y(\theta) = \bar{a}(1 - \cos\theta) + (vt - \bar{a}\theta)\sin\theta, \quad 0 \leq \theta \leq vt/\bar{a} \quad (4)$$

Unlike C-paths, where the reachability region is spanned by curves of all turn radii $a \geq \bar{a}$, [9] showed that the optimal turn radius for maximum reach of CS-paths is always the minimum turn radius $\bar{a}$. While [9] offered no derivation for (4), the following is how one can derive the equation.

Consider the vehicle to make a turn of angle $\theta$ starting at $O$ and ending at $G$, then continuing in a straight line until reaching $F$, as depicted in Fig. 4. The location of $F$ is simply:

$$F = O + \overrightarrow{OC} + \overrightarrow{CG} + \overrightarrow{GF} \quad (5)$$

This equation can be broken into two pieces:

$$F = O + \underbrace{\begin{pmatrix} 0 \\ \bar{a} \end{pmatrix} + \begin{pmatrix} \bar{a}\sin\theta \\ -\bar{a}\cos\theta \end{pmatrix}}_{\text{turn}} + \underbrace{(vt - \bar{a}\theta)\begin{pmatrix} \cos\theta \\ \sin\theta \end{pmatrix}}_{\text{straight}} \quad (6)$$

The horizontal and vertical lines at the point $G$ in Fig. 4 are drawn to assist the reader in this explanation. The angle $\theta$ is bound to be in the domain: $0 \leq \theta \leq vt/\bar{a}$. Sweeping over these angles, the top-half of the CS-boundary can be generated. The bottom half of the curve is obtained by reflecting (4) over the x-axis, and the complete boundary is shown in Fig. 5 with a red line highlighting a single trajectory. Alternatively, the boundary can be defined in a single parametric expression with $-vt/\bar{a} \leq \theta \leq vt/\bar{a}$ by modifying (4) to:

$$x(\theta) = \bar{a}\sin|\theta| + (vt - \bar{a}|\theta|)\cos\theta$$
$$y(\theta) = \text{sign}(\theta)(\bar{a}(1 - \cos\theta) + (vt - \bar{a}|\theta|)\sin|\theta|) \quad (7)$$

As before, the circles were added to the figure separately, and, similar to C-paths, the locus travels inside its outer boundary when $vt$ exceeds $\left(1 + \frac{3\pi}{2}\right)\bar{a}$. The reachability region is all the space contained radially between this bounding curve and the minimum turn radius circles until this criteria is reached.

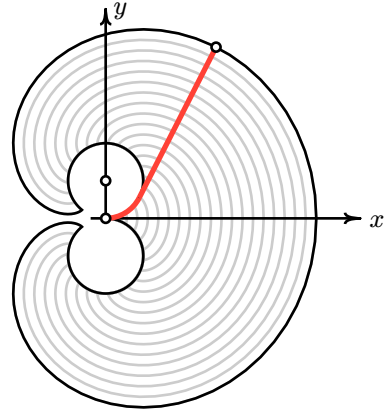

Fig. 5. A vehicle (red) takes a turn-straight trajectory for a fixed amount of time. The locus of all possible terminal conditions is shown in dark black, in grey are iso-cost lines for other time intervals to illustrates how the frontier grows with time.

Unlike with C-paths, where reaching the inside of the minimum turn circles is impossible, for CS-paths, the inside of the circles is first reachable when $vt \geq \left(-1 + \sqrt{3} + \frac{5\pi}{3}\right)\bar{a}$, and the entirety of the circles is reachable when $vt \geq \left(1 + \sqrt{3} + \frac{5\pi}{3}\right)\bar{a}$. There is some subtlety in determining which space internal to them is reachable if $vt$ is between these values, though this is likely unnecessary, as non-pathological navigation problems will not involve this region at all.

## C. Comparison of Paths

If a vehicle were holonomic, it could change its instantaneous heading at any time step. The reachability for such a vehicle with limited range and fixed speed would be a circle centered at the origin with radius $vt$. Both (3) and (4) reduce to this as $\bar{a}$ goes to zero. In Fig. 6, the two reachability regions previously described are compared against that of the reachability region of a holonomic vehicle.

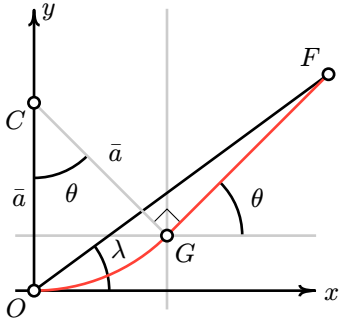

Fig. 4. A vehicle initially pointing along the $x$-axis takes a turn-straight trajectory, turning until reaching a point $T$ and goes straight/terminates at a point $F$.

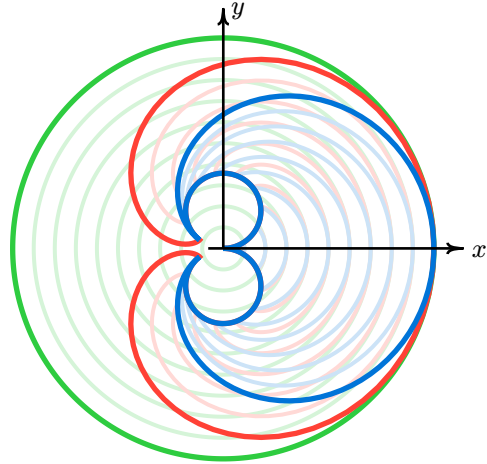

Fig. 6. A comparison of three different reachability regions for a fixed amount of time and speed is shown. In green is the reachability region for a holonomic vehicle. In red is the reachability region for a vehicle taking a turn-straight path. In blue is the reachability region for a turn-constrained vehicle making constant turn.

## III. ENGAGEMENT ZONES

This section reformulates the previously described reachability regions into engagement zones for use in path planning problems. A key feature of an EZ, detailed in [3], is that it is the same size and shape as a pursuer's reachability region, but translated in space based on the target's speed and heading. If the target alters its course, the EZ changes dynamically, rotating about the pursuer's position by the same angle as the heading change of the target. If the target alters its speed, the EZ translates along the heading vector of the target, either toward or away from the pursuer depending on whether the target is slowing down or speeding up.

More specifically, a pursuer's EZ results from its reachability region being translated a distance $\mu R$ in the direction opposite the target's heading, where $\mu$ is the speed ratio $v_T/v_P$ [3]. If the target changes its heading while maintaining speed, the EZ will "orbit" about the pursuer's location such that the translation is always opposite the target's new heading. This occurs without changing the EZ's own relative orientation, which is determined solely by the heading of the pursuer. A speed change on the part of the target changes the magnitude of the translational shift through altering $\mu$.

### A. Turn-Only Basic Engagement Zone (CBEZ)

If a pursuing vehicle is turn constrained, has limited range, and only executes a single turn (subject to the turn constraint), the reachable region for that vehicle follows (3). The reachability region does not account for the time of flight of the pursuer nor does it account for a target's current speed and course. By backward propagating a target's trajectory from the reachability region for the time that the pursuer moves, the engagement zone for this turn-only case can be obtained. The EZ and this procedure are graphically shown in Fig. 7 and detailed below.

As (3) assumes an initial heading in the x-direction, it must be rotated by the heading of the pursuer relative to the x-axis. Since (3) is already in polar form, this rotation can be easily represented by the shift:

$$r(\lambda) = \frac{vt}{\lambda - \psi_P} \sin(\lambda - \psi_P) \qquad (8)$$
$$-\frac{vt}{2\bar{a}} + \psi_P \leq \lambda \leq \frac{vt}{2\bar{a}} + \psi_P$$

Then the resulting shape is translated a distance $\mu R$ opposite the target's heading. For the purposes of path planning, a constraint must be implemented based on this EZ border. The following is how one can be defined.

As in Fig. 7, let $P'$ be the pursuer's position $P$ shifted based on the target's heading $\psi_T$ relative to the x-axis:

$$\overrightarrow{PP'} = \begin{pmatrix} -\mu R \cos \psi_T \\ -\mu R \sin \psi_T \end{pmatrix} \qquad (9)$$

Note that $\psi_T = 0$ in Fig. 7 for simplicity. Track the target's position $T$ relative to this new origin. The target will be at some polar angle $\lambda'$ relative to $P'$. If $T$ lies on the border of the EZ:

$$|\overrightarrow{P'T}| = \frac{vt}{\lambda' - \psi_P} \sin(\lambda' - \psi_P) \qquad (10)$$

Where the pursuer heading $\psi_P$ is accounting for a rotation of the EZ. Generally, define $|\overrightarrow{P'T}| \equiv d'$. To stay out of the EZ, the target must achieve:

$$d' \geq \frac{vt}{\lambda' - \psi_P} \sin(\lambda' - \psi_P) \qquad (11)$$

With consideration to the bounds $-\pi \leq (\lambda' - \psi_P) \leq \pi$. Some care must be taken to recast an arbitrary $\lambda' - \psi_P$ within these bounds, as it is not correct to simply clamp the value, rather, it must be converted to an equivalent unit circle angle between $-\pi$ and $\pi$.

### B. Turn-Straight Basic Engagement Zone (CSBEZ)

Formulating the EZ for turn-straight motion is more difficult, as the reachability region is not easily cast in polar form due to the complex relationship between $\theta$ and $\lambda$. This makes defining a constraint that keeps a vehicle out of the EZ, akin to (11), challenging. A semi-analytic approach is outlined below.

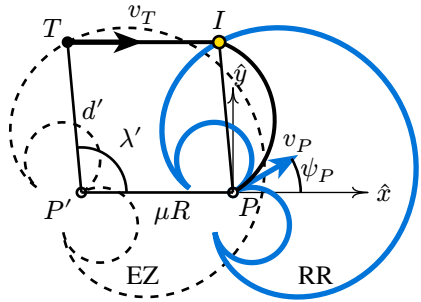

Fig. 7. A target moves to the right while a pursuer with limited range and turn rate is currently aiming upward and to the right. The pursuer can only make soft constant turns. The RR of the pursuer is in blue. The target is currently on the EZ of the pursuer shown in dashed black.

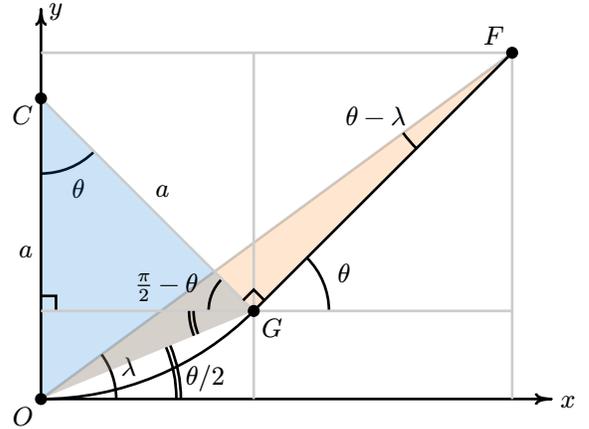

Fig. 8. CS-path detail for conversion to polar coordinates

First, it is necessary to derive a relationship between the turn angle $\theta$ and the polar angle $\lambda$ of a given trajectory's end point. Consider the $\triangle OFG$ in Fig. 8. Using the Law of Cosines, the following can be said about the triangle:

$$\overline{OG}^2 = \overline{OF}^2 + \overline{GF}^2 - 2\overline{OF}\,\overline{GF}\cos(\theta - \lambda) \quad (12)$$

Solving (12) for $\theta$:

$$\theta = \lambda + \cos^{-1}\left(\frac{\overline{OF}^2 + \overline{GF}^2 - \overline{OG}^2}{2\overline{OF}\,\overline{GF}}\right) \quad (13)$$

Where:
$$\begin{aligned}\overline{OF} &= |(a\sin\theta + (vt - a\theta)\cos\theta)\hat{x} + \\ &\quad (a(1-\cos\theta) + (vt - a\theta)\sin\theta)\hat{y}| \\ \overline{GF} &= vt - a\theta \\ \overline{OG} &= 2a\sin(\theta/2)\end{aligned} \quad (14)$$

As with the turn-only EZ, we would like to find the radial distance from the center $P'$ of the turn-straight EZ to its boundary, along the line from $P'$ to the target's position T, as this represents the minimum allowable distance to stay outside the EZ (shown in Fig. 9). This is given by $\overline{OF}$ in (14), which requires knowledge of $\theta$ to evaluate. It is conceivable to use a basic root-finder on (13) to find $\theta$ with an initial guess of $\lambda$. This is straightforward as $\cos^{-1}(.)$ is a monotonically increasing function with a range of $[0, \pi]$, and thus we must have $\lambda \leq \theta \leq \lambda + \pi$.

Importantly, this must be done with the polar angle $\lambda'$ of the target's position relative to the center of the EZ as described in the previous section. To account for a nonzero pursuer heading $\psi_P$, the same procedure should instead be done with $\lambda' - \psi_P$, converted to an equivalent unit circle angle between $-\pi$ and $\pi$. Due to the significant added complexity involved with the turn-straight EZ, the path planning example in this paper only uses the turn-only formulation.

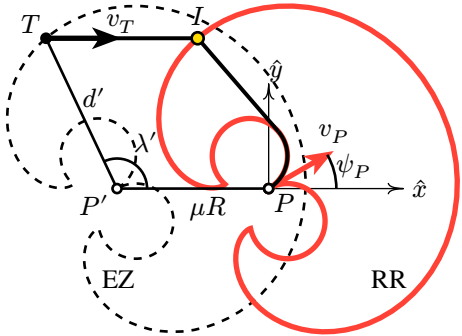

Fig. 9. A target moves to the right while a pursuer with limited range and turn rate is currently aiming upward and to the right. The pursuer can take turns and straight-line trajectories. The RR of the pursuer is in red. The target is currently on the EZ of the turn-constrained pursuer shown in dashed black.

## IV. APPLICATION

In this section, the turn-only basic engagement zone (CBEZ) is utilized within a path planning problem. Here, a vehicle moving at constant speed with simple motion (i.e., single integrator kinematics) starts from a specified location and wishes to reach a target location in minimum time while avoiding entering the CBEZ of a pursuer located at the origin. The CBEZ arises due to the presence of a threat which is assumed to move with speed greater than that of the vehicle and is subject to a turn-rate constraint. Over the course of the path, it is assumed that the threat never launches. Nevertheless, the vehicle seeks to avoid entering the CBEZ lest the threat have an opportunity to launch and intercept it. Let $\psi_T$ represent the vehicle's instantaneous heading, over which it has control. This problem is formulated as follows:

$$\begin{aligned}&\min_{\psi_{T(t)}} t_f \\ \text{s.t.} \quad &\dot{x} = \mu v \cos(\psi_T), \quad \dot{y} = \mu v \sin(\psi_T) \\ &x(0) = x_0, \quad y(0) = y_0 \\ &x(t_f) = x_f, \quad y(t_f) = y_f \\ &d' \geq \frac{vt}{\lambda' - \psi_P}\sin(\lambda' - \psi_P)\end{aligned} \quad (15)$$

where $d'$ and $\lambda'$ are defined in (11) and (10), respectively, and $x_0, y_0, x_f,$ and $y_f$ are the specified initial and final locations.

The solution to (15) is computed numerically using the the even collocation method (c.f., e.g., [23]). Additionally, the solution is compared against the solution to an augmented version of (15) wherein the EZ constraint is replaced by the BEZ constraint from [3]. In the augmented version, the pursuing vehicle is assumed to have no turn-rate constraint and can freely move in any direction at any moment. Finally, both of these solutions are compared against a nominal trajectory based on circumnavigating a circle of radius $R$. For these results, the parameters were selected as $\mu = 0.9$, $\bar{a} = 0.25$, $v = 1$, $\psi_P = \pi$, and $R = t = \frac{\pi}{2}$.

Table I contains the results for the different trajectories. The trajectories themselves are shown in Fig. 10 along with the reachable regions for the simple motion and turn-constrained pursuer vehicles. Fig. 11 shows the CBEZ-avoiding trajectory at various points in time in order to show the CBEZ and how it's position changes w.r.t. the target vehicle's heading.

Compared to circumnavigating the circle of radius $R$, the BEZ-based path plan yields a 2.24% time savings, while the CBEZ-based path plan yields a 5.19% savings. The CBEZ path is 3.02% shorter than the BEZ path. This highlights the fact that the assumption of the threat having simple motion may be overconservative for the case in which the threat is actually turn-constrained. Put differently, knowledge of the threat's capabilities and motion behavior can be exploited by the vehicle in order to plan shorter paths while still maintaining the same level of risk.

TABLE I. TRAVEL TIMES FOR THE TRAJECTORIES

| Trajectory | Time | % Improvement |
|---|---|---|
| Nominal | 9.803 | - |
| BEZ | 9.584 | 2.24% |
| CBEZ | 9.294 | 5.19% |

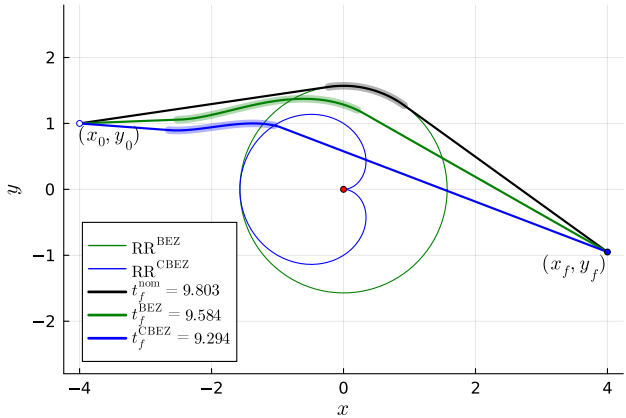

Fig. 10. Example path plans which start at $(x_0, y_0)$ and go to $(x_f, y_f)$ in minimum time while avoiding entering the EZ. The black path is the nominal trajectory which circumnavigates the circle of radius $R$. The green path avoids the simple-motion BEZ (as in [3]) while the blue path avoids the CBEZ. Highlighted sections indicate where the associated constraint is active.

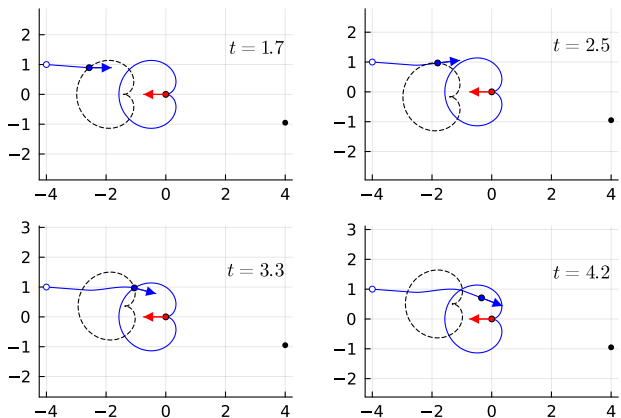

Fig. 11. Snapshots of the CBEZ-avoiding trajectory from Fig. 10 with the instantaneous CBEZ shown in black.

## V. CONCLUSIONS

This work formulates two novel engagement zones for turn-constrained vehicles: CBEZ and CSBEZ. The CBEZ corresponds to turn-only paths and CSBEZ corresponds to turn-straight paths of a Dubins' vehicle. The models, derived from first principles begin with reachability analysis followed by derivations of corresponding basic engagement zones. A simulation using nonlinear programming demonstrates utility for path planning around turn constrained pursuers with limited range. Moreover, by considering the turn constraint of the pursuer, the vehicle is able to reduce its overall path length. Future work considers utilizing these models for engagements involving many pursuers and or evaders. The semi-analytic nature of the CSBEZ may be challenging and future methods such as sampling methods may prove useful for this model.